\documentclass[12pt,reqno]{amsart}
\usepackage[margin=1in]{geometry}
\usepackage[utf8]{inputenc}
\usepackage{setspace}
\usepackage{graphicx}
\usepackage{xcolor}
\usepackage{amsmath, mathtools}
\usepackage{amsthm, amssymb}
\usepackage{enumitem}
\usepackage{caption}
\usepackage{subcaption}

\usepackage{todonotes}
\usepackage{tikz}
\usepackage[normalem]{ulem}

\newtheorem{theorem}{Theorem}

\newtheorem{lemma}{Lemma}

\usepackage{url}
\usepackage{hyperref}
\hypersetup{
    colorlinks,
    linkcolor={red!60!black},
    citecolor={green!60!black},
    urlcolor={blue!60!black}
}

\tikzstyle{vert}=[shape=circle,draw=black,fill=black, inner sep=.75mm]
\tikzstyle{fixed}=[shape=rectangle,draw=black,fill=white, inner sep=1.2mm]
\tikzstyle{uncolored}=[dashed,thick]
\tikzstyle{uncolored2}=[dotted,thick]
\tikzstyle{red?}=[dashed,thick,color=red]
\tikzstyle{blue?}=[dashed,thick,color=blue]
\tikzstyle{purple?}=[dashed,thick,color=purple]
\tikzstyle{purple}=[solid,thick,color=purple]
\tikzstyle{red}=[solid,thick,color=red]
\tikzstyle{blue}=[solid,thick,color=blue]
\tikzstyle{green}=[solid,thick,color=green]
\tikzstyle{green?}=[dashed,thick,color=green]

\newcommand{\prob}{\mathbb{P}}
\newcommand{\E}{\mathbb{E}}

\newcommand{\rbrac}[1]{\left(#1\right)} %round brackets
\DeclareMathOperator{\supp}{supp}
\DeclareMathOperator{\wt}{wt}
\DeclareMathOperator{\di}{d}
\newcommand{\mc}[1]{\mathcal{#1}}

\renewcommand{\l}{\ell}

\newcommand{\vx}{{\vec{x}}}
\newcommand{\vy}{{\vec{y}}}
\newcommand{\vw}{{\vec{w}}}

\allowdisplaybreaks

\title{Asymptotically optimal constant weight codes with even distance}

\author{Patrick Bennett}
\address{Department of Mathematics, Western Michigan University, Kalamazoo, MI, USA}
\thanks{The author was supported in part by Simons Foundation Grant \#426894.}
\email{\tt patrick.bennett@wmich.edu}

\begin{document}

\maketitle

\begin{abstract}
    A $q$-ary code $C$ of length $n$ is a set of $n$-dimensional vectors (code words) with entries in $\{0, \ldots, q-1\}$. We say $C$ has constant weight $w$ if each code word has exactly $w$ nonzero entries. We say $C$ has minimum distance $d$ if any two distinct code words in $C$ differ in at least $d$ entries. We let $A_q(n, d, w)$ be the largest possible cardinality of any $q$-ary code of length $n$ with constant weight $w$ and minimum distance $d$. Very recently, Liu and Shangguan gave an asymptotically sharp estimate for $A_q(n, d, w)$ where $q, d, w$ are fixed, $d$ is odd and $n \rightarrow \infty$. In this note we answer a question of Liu and Shangguan by obtaining such an estimate in the case where $d$ is even. 
\end{abstract}

\section{Introduction}

Let $\Sigma_q:=\{0, 1, \ldots, q-1\}$, and let $\Sigma_q^n$ be the set of $n$-dimensional vectors with entries in $\Sigma_q$. We call the elements of $\Sigma_q^n$  {\em code words}. For $\vx \in \Sigma_q^n$, we let $x_i$ be the $i^{th}$ entry of $\vx$ and we define the {\em support} and {\em weight} of $\vx$ as
\begin{equation}
  \supp(\vx):= \{i: x_i \neq 0\}, \qquad   \wt(\vx):= |\supp(\vx)|.
\end{equation}
For $\vx, \vy \in \Sigma_q^n$, the {\em Hamming distance} between $\vx$ and $\vy$ is  
\begin{equation}
    \di(\vx, \vy):= \big|\{i: x_i \neq y_i\}\big|.
\end{equation}
We say $C$ is a {\em constant weight code with weight $w$} if $\wt(\vx)=w$ for all $\vx \in C$. We say $C$ has {\em minimum distance $d$} if $\di(\vx, \vy) \ge d$ for all distinct $\vx, \vy \in C$. We let $A_q(n, d, w)$ be the largest possible cardinality of any $q$-ary code of length $n$ with constant weight $w$ and minimum distance $d$.

Constant weight codes have several applications (see \cite{7, 14, CF07}) and are well-studied (see \cite{17, 19, 10, 2, 13, 18, 5, 12, 20,  15, 3, 16, 11, 9,  6, 8}). Much of the work centers on estimating $A_q(n, d, w)$. The following upper bound appears, for example, in Chee and Ling \cite{13}:
\begin{theorem}\label{thm:johnson}
    Let $t=\lceil \frac{2w-d+1}{2} \rceil$. Then
    \begin{equation*}
        A_q(n,d,w)\le \left\{ \begin{array}{cl}
                   \frac{(q-1)^t \binom{n}{t}}{\binom{w}{t}},& \text{if $d$ is odd;} \\
                    \frac{(q-1)^{t-1} \binom{n}{t}}{\binom{w}{t}}, & \text{if $d$ is even}.
                \end{array}\right.
    \end{equation*}
\end{theorem}
The above bounds follow from iteratively applying some classical bounds due to Johnson \cite{J62}. However they also follow from the Pigeonhole Principle (we will see in Section \ref{sec:motivation}). Very recently, Liu and Shangguan \cite{LS25} proved that when $d$ is odd, Theorem \ref{thm:johnson} is asymptotically optimal, and asked whether the same holds when $d$ is even. In this note we will answer that question in the affirmative:
\begin{theorem}\label{thm:main}
Fix $q \ge 2, w \ge 1$ and an even number $d$,  $2 \le d \le 2w$. Let $t:= \lceil \frac{2w-d+1}{2} \rceil = w - \frac12 d + 1$. Then
  \begin{equation}
    A_q(n, d, w) = (1+o(1))  \frac{(q-1)^{t-1} \binom{n}{t}}{\binom{w}{t}}.
\end{equation}  
\end{theorem}
 In Section \ref{sec:tools} we introduce some basic terminology and tools. In Section \ref{sec:motivation} we discuss the idea of the proof of Theorem \ref{thm:main} at a high level. In Section \ref{sec:proof} we prove Theorem \ref{thm:main}. 

\section{Terminology and Tools}\label{sec:tools}

We use standard asymptotic notation $O(-), o(-), \Omega(-)$. We also say $f(n) \sim g(n)$ if $f(n)=(1+o(1))g(n)$. All asymptotics are as $n \rightarrow \infty$ while $q, d, w$ are all fixed. For a positive integer $m$ we let $[m]=\{1, \ldots, m\}$. We say a set of size $t$ is a $t$-set, and for a set $S$ of size at least $t$, we write $\binom St$ for the collection of all $t$-sets contained in $S$. 

A {\em hypergraph} $\mc{H}$ on a vertex set $V=V(\mc{H})$ is a collection of subsets of $V$. The elements of $\mc{H}$ are called {\em edges}. We say $\mc{H}$ is $\l$-bounded if every edge has size at most $\l$. A {\it matching} in $\mc{H}$ is a set of pairwise disjoint edges. The {\it matching number} of $\mc{H}$, denoted by $\nu(\mc{H})$, is the maximum possible size of a matching in $\mc{H}$. We say that a function $f:E(\mathcal{H})\rightarrow [0,1]$ is a {\it fractional matching} for $\mathcal{H}$ if for every $v\in V(\mathcal{H})$, 
$$\sum_{e\in \mathcal{H}:~v\in e}f(e)\le 1.$$ 
For a fractional matching $f$ for $\mathcal{H}$, let $f(\mathcal{H}):=\sum_{e\in E(\mathcal{H})}f(e)$ and
$$\alpha(f):=\max_{x\neq y \in V(\mathcal{H})} \sum_{e\in E(\mathcal{H}):~x,y\in e} f(e).$$

Kahn \cite{K96} showed that under certain conditions, the existence of a ``large'' fractional matching implies the existence of a ``large'' matching:
\begin{lemma}[Theorem 1.2 in \cite{K96}]\label{lem:Kahn}
    For every integer $\l\ge 1$ and every real $\epsilon>0$, there exists a real $\sigma>0$ such that whenever $\mc{H}$ is a $\l$-bounded hypergraph and $f$ is a fractional matching of $\mc{H}$ with $\alpha(f)<\sigma$, then $\nu(\mc{H})>(1-\epsilon)f(\mathcal{H}).$
\end{lemma}
Our proof of Theorem \ref{thm:main} will involve applying Lemma \ref{lem:Kahn} to a carefully constructed hypergraph $\mc{H}$. Part of the construction of $\mc{H}$ will be random, and we will show that this randomly constructed $\mc{H}$ satisfies certain properties {\em with high probability (w.h.p.)}, i.e. with probability tending to 1 as $n \rightarrow \infty$. 
To obtain our probability bounds, we will use the following concentration inequality due to McDiarmid~\cite{M89}: 
\begin{lemma}\label{thm:McDiarmid}
Suppose $B_1,\ldots,B_m$ are independent random variables. Suppose $B$ is a real-valued random variable determined by $B_1,\ldots, B_m$ such that changing the outcome of $B_i$ changes $B$ by at most $b_i$ for all $i\in [m]$. Then, for all $\lambda \geq 0$, we have \[\prob[|B-\E[B]|\geq \lambda ]\leq 2\exp\left(-\frac{2\lambda^2}{\sum_{i=1}^m b_i^2}\right).\]
\end{lemma}

\section{Motivation}\label{sec:motivation}

First, let us see a proof of the upper bound which illuminates what one must do to construct a code nearly that large. We start with the case of odd $d$. Odd $d$ is the case where Liu and Shangguan \cite{LS25} recently proved Theorem \ref{thm:main} is asymptotically optimal, and we will briefly sketch a (technically slightly different) proof of that result. Then we will move on to the case of even $d$ and discuss the complication that is not present for odd $d$. 
\begin{proof}[Proof of Theorem \ref{thm:johnson}, odd $d$] 
Let $C$ be a $q$-ary constant weight code of weight $w$ and minimum distance $d$. We use the Pigeonhole Principle. For each $\vx \in C$, there are $\binom wt$ many $t$-sets in $\supp(\vx)$, and with each such $t$-set we associate a vector in $[q-1]^t$ specifying the entries of $\vx$ on $\supp(\vx)$. If some $\vy \neq \vx$ in $C$ had $|\supp(\vy) \cap \supp(\vx)| \ge t$ and $\vy, \vx$ both had the same vector of nonzero entries on some $t$-set, then we would have $\di(\vx, \vy) \le 2(w-t) < d$, a contradiction. The total number of choices for a $t$-set in $[n]$ and a  a vector in $[q-1]^t$ is $(q-1)^t \binom nt $. Thus we have 
\[
|C| \le   \frac{(q-1)^t \binom{n}{t}}{\binom{w}{t}},
\]
so we are done in this case. 
\end{proof}
Roughly speaking, for the Pigeonhole Principle to be (nearly) tight, we need to put one pigeon in (almost) every hole. In other words, for the above proof to give a nearly tight bound, it must be the case that for almost every $t$-set of indices $T$ in $[n]$ and almost every possible vector of entries $(q-1)^t$, there is some $\vx \in C$ having those specified entries on the set of indices $T$. This is the idea behind Liu and Shangguan's \cite{LS25} proof, which is almost the same as the proof we now sketch. 
\begin{proof}[Sketch of Liu and Shangguan's result \cite{LS25}]
 Consider the following hypergraph $\mc{H}$. Let $V_1 :=\binom{[n]}{t+1}$. Let $V_2$ be the family of all sets of the form $\{(i_1, x_{i_1}), \ldots (i_t, x_{i_t})\}$ such that the $i_j \in [n]$ are all distinct and the $x_{i_j}$ are in $[q-1]$. The vertex set of our hypergraph will be $V(\mc{H}):= V_1 \cup V_2$. 
 
 We now describe the edges of $\mc{H}$. Let $X$ be the set of all code words  $\vx =(x_1, \ldots x_n) \in \Sigma_q^n$ having weight $w$. For each $\vx \in X$ we will have a hyperedge $e_\vx \in \mc{H}$. Specifically, we will have $e_\vx \cap V_1 = \binom{\supp(\vx)}{t+1}$. Now supposing that $\supp(\vx) = \{i_1, \ldots, i_w\}$, we will have $e_\vx \cap V_2 = \binom{\{(i_1, x_{i_1}), \ldots, (i_w, x_{i_w})\}}{t}$.

 Note that there is a one-to-one correspondence between matchings $M \subseteq \mc{H}$ and $q$-ary codes $C$ with constant weight $w$ and minimum distance $d$. Indeed, the vertices in $V_1$ enforce that for any two code words $\vx, \vy \in C$, we have $|\supp(\vx) \cap \supp(\vy)| \le t$, and the vertices in $V_2$ enforce that whenever we do have $|\supp(\vx) \cap \supp(\vy)| = t$, the code words $\vx, \vy$ do not completely agree on their common support. 

 An application of Lemma \ref{lem:Kahn} to $\mc{H}$ (details omitted) yields a matching that covers almost all of $V_2$, which then corresponds to a code proving Liu and Shangguan's lower bound \cite{LS25}. 
    
\end{proof}

We now turn to the case of even $d$. First we show how the upper bound follows from the Pigeonhole Principle. 

\begin{proof}[Proof of Theorem \ref{thm:johnson}, even $d$] 
For each $\vx \in C$, there are $\binom wt$ many $t$-sets $T \subseteq \supp(\vx)$, and with each such $T$ we associate a vector in $[q-1]^{t-1}$ specifying the entries of $\vx$ on {\em some arbitrary $t-1$ indices} in $T$. For definiteness, we could say we always take the lowest $t-1$ indices of $T$. If some $\vy \neq \vx$ in $C$ had $|T|=t$ where $T=\supp(\vy) \cap \supp(\vx)$ and $\vy, \vx$ both had the same vector of nonzero entries on the associated $(t-1)$-set, then we would have $\di(\vx, \vy) \le 2(w-t)+1 < d$, a contradiction. Thus we have 
\[
|C| \le   \frac{(q-1)^{t-1} \binom{n}{t}}{\binom{w}{t}}.
\]
\end{proof}
Of course it seems a bit strange that in the proof above, we only use a vector in $[q-1]^{t-1}$ obtained by arbitrarily discarding information. For the result to be nearly tight, it must be nearly tight for every possible way of making this arbitrary choice. In other words, for almost every $t$-set of indices $T$ and $(t-1)$-subset $T'\subseteq T$, we see every possible vector in $[q-1]^{t-1}$ on $T'$. We also have to guarantee that whenever we have two code words $\vx, \vy$ with $\supp(\vx) \cap \supp(\vy) = T$, we have that $\vx, \vy$ disagree in two entries (otherwise $\di(\vx, \vy) <d$). To accomplish this, for each set $T$ we will impose an equation that $\vx$ must satisfy if $T \subseteq \supp(\vx)$. Specifically, for each $t$-set of indices $T$ and $\vx$ containing $T$ in its support, we will require 
 \begin{equation}\label{eqn:trick}
        \sum_{i \in T} x_i = B_{T} \pmod{q-1}
    \end{equation}
for some appropriately chosen $B_T$. First, note that for every $(t-1)$-subset $T'\subseteq T$ it is possible to get all $(q-1)^{t-1}$ vectors of nonzero entries on $T'$ since we can set these $t-1$ entries however we like and choose the remaining entry (on the single index in $T\setminus T'$) to satisfy \eqref{eqn:trick}. Next, note that if $\vx, \vy$ both contain $T$ in their support, then it is not possible for $\vx, \vy$ to disagree on exactly one entry in $T$. Indeed, by \eqref{eqn:trick} we would have $\sum_{i \in T} x_i = \sum_{i \in T} y_i$ which is not possible when $x_i=y_i$ for all but one index $i$.

At this point we have a plan. We will apply Lemma \ref{lem:Kahn} to a hypergraph $\mc{H}$ which we will construct so that a matching in $\mc{H}$ corresponds to a set of code words that satisfy the equations \eqref{eqn:trick}. However the choice of the values $B_T$ must be done so that $\mc{H}$ will be have nice properties for Lemma \ref{lem:Kahn}. Choosing, for example, $B_T=0$ for all $T$ would be a bad choice. As we will see, a random choice for the $B_T$ works well.  

\section{Proof of Theorem \ref{thm:main}}\label{sec:proof}

\begin{proof}
    For each $T \in \binom{[n]}{t}$ independently, we choose $B_T \in [q-1]$ uniformly at random. 
 Let $V_1 :=\binom{[n]}{t+1}$. Let $V_2$ be the family of all sets of the form $\{(i_1, x_{i_1}), \ldots (i_t, x_{i_t})\}$ such that the $i_j \in [n]$ are all distinct and where the $x_{i_j} \in [q-1]$ satisfy the equation  \eqref{eqn:trick} where $T=\{i_1, \ldots, i_t\}$. Our vertex set will be $V(\mc{H})=V_1 \cup V_2$.

     We now describe the edges of $\mc{H}$. Let $X$ be the set of all code words  $\vx =(x_1, \ldots x_n) \in [q-1]^n$ having weight $w$ and such that for all $T \in \binom{\supp(x)}{t}$ we have \eqref{eqn:trick}. For each $\vx \in X$ we will have a hyperedge $e_\vx \in \mc{H}$. Specifically, $e_\vx \cap V_1 = \binom{\supp(\vx)}{t+1}$. Now supposing that $\supp(\vx) = \{i_1, \ldots, i_w\}$, we will have $e_\vx \cap V_2 = \binom{\{(i_1, x_{i_1}), \ldots, (i_w, x_{i_w})\}}{t}$. So each edge has $\binom{w}{t+1}$ vertices from $V_1$ and $\binom{w}{t}$ vertices from $V_2$. 

     We estimate the degrees of vertices in $\mc{H}$. A vertex in $V_1$ has degree $O(n^{w-t-1})$. We claim that w.h.p. (with respect to the randomness in $\mc{H}$), the degree of every vertex in $V_2$ is $\sim (q-1)^{w-t-\binom{w}{t}+1} \binom{n-t}{w-t}$. Indeed, we consider a fixed tuple $v=\{(i_1, x_{i_1}), \ldots (i_t, x_{i_t})\}$ (which may or may not be in $V_2$ since we do not assume that \eqref{eqn:trick} holds). We estimate the number $D$ of code words $\vx$ such that $e_{\vx} \ni v$ satisfies all the requirements to be in $\mc{H}$ (except for possibly the requirement that $v \in V_2$, i.e.~an equation of the form \eqref{eqn:trick}). We have 
     \[
     \E[D] =  (q-1)^{w-t-\binom{w}{t}+1} \binom{n-t}{w-t}.
     \]
     Indeed, to complete $\vx$ given $v$ we have to choose $w-t$ more elements of $[n]$ to be in $\supp(\vx)$, then we have $(q-1)^{w-t}$ choices for the entries of $\vx$ there, and finally for each of the $\binom{w}{t}-1$ additional equations of the form \eqref{eqn:trick} which $\vx$ has to satisfy, there is a probability of $(q-1)^{-1}$ that the equation is satisfied. Note that $D$ is determined by the random variables $B_{S'}, S' \neq \{i_1, \ldots, i_t\}$. Consider some $S'$ such that $|S' \setminus \{i_1, \ldots, i_t\}|=k \ge 1$ and note that the number of such $S'$ is $O(n^k)$. Then changing the value $B_{S'}$ could change $D$ by at most $O(n^{w-t-k})$. Indeed, the only vectors $\vx$ that could be affected here must have $\supp(\vx) \supseteq \{i_1, \ldots, i_t\} \cup S'$, which has size $t+k$. We apply Lemma \ref{thm:McDiarmid} with $\lambda := n^{w-t-1/3} = o(\E[D])$ to obtain
     \[\prob[|D-\E[D]|\geq \lambda ]\leq \exp\left(-\Omega \rbrac{\frac{n^{2w-2t-2/3}}{\sum_{k=1}^t n^k (n^{w-t-k})^2}}\right) = \exp\rbrac{-\Omega\rbrac{n^{1/3}}}.\]
     Now by the union bound over all choices for $v\in V_2$, the probability that any such $v$ has a degree outside the interval $\E[D] \pm \lambda$ is at most $n^{O(1)} \exp\rbrac{-\Omega\rbrac{n^{1/3}}}= o(1).$ We henceforth assume that $\mc{H}$ is chosen to be one of the more likely outcomes of the random experiment, i.e. that every vertex in $V_2$ has degree $\sim(q-1)^{w-t-\binom{w}{t}+1} \binom{n-t}{w-t} $.

     We will now apply Lemma \ref{lem:Kahn}. Let $f:\mc{H}\rightarrow [0, 1]$ be the constant function $f(e) = \phi := (\Delta(\mc{H}))^{-1}  \sim\rbrac{(q-1)^{w-t-\binom{w}{t}+1} \binom{n-t}{w-t}}^{-1} $, so trivially $f$ is a fractional matching. Since any two vertices of $\mc{H}$ are in $O(n^{w-t-1})$ edges together, we have $\alpha(f) = O(n^{w-t-1} \cdot \phi ) = O(n^{-1})$. Note that $|V_2| =(q-1)^{t-1}\binom{n}{t}$ and recall that every edge of $\mc{H}$ has $\binom wt$ vertices in $V_2$. Thus $\mc{H}$ has a matching $M$ of size asymptotically 
     \[
      f(\mc{H}) =  |\mc{H}| \phi = \frac{\sum_{v \in V_2} deg(v)}{\binom wt \Delta(\mc{H})}\sim  \frac{(q-1)^{t-1}\binom{n}{t}}{\binom{w}{t}}.
     \]
     By construction, this matching $M$ corresponds to a code $C$ of the same size. Indeed, consider two code words $\vx, \vy$. We cannot have $|\supp(\vx)\cap \supp(\vy)|>t$, or else the edges $e_{\vx}, e_{\vy}$ would intersect in $V_1$. If $|\supp(\vx)\cap \supp(\vy)|=t$ then $\vx, \vy$ must both satisfy \eqref{eqn:trick} where $T=\supp(\vx)\cap \supp(\vy)$. In this case $\vx, \vy$ also cannot completely agree on $T$ or else they would intersect in $V_2$. Therefore $\vx, \vy$ must actually disagree in at least two entries of $T$ since they both satisfy \eqref{eqn:trick}, and so $\di(\vx, \vy) \ge 2(w-t)+2 = d$. Finally, if $|\supp(\vx)\cap \supp(\vy)|\le t-1$ then $\di(\vx, \vy) \ge 2(w-t+1) = d$. This completes the proof. 
\end{proof}

\section{Concluding remarks}

In this note we answered a question of Liu and Shangguan \cite{LS25} by constructing a constant weight code that asymptotically matches the Johnson bound \cite{J62}. In their paper \cite{LS25}, Liu and Shangguan also studied another type of code called a {\em constant composition code}. For this type of code there is some fixed vector $\vw\in \mathbb{N}^{q-1}$ such that for all code words $\vx$ we have $|\{i: x_i=j\}|=w_j$ for all $1 \le j \le q-1$. Liu and Shangguan \cite{LS25} showed that the Johnson bound for constant composition codes is tight for odd $d$, and asked whether it also holds for even $d$. We tried to answer that question by adapting our proof to that setting, but we could not seem to asymptotically match the Johnson bound.

\bibliographystyle{abbrv}

\begin{thebibliography}{10}

\bibitem{17}
E.~Agrell, A.~Vardy, and K.~Zeger.
\newblock Upper bounds for constant-weight codes.
\newblock {\em IEEE Trans. Inform. Theory}, 46(7):2373--2395, 2000.

\bibitem{19}
A.~E. Brouwer, J.~B. Shearer, N.~J.~A. Sloane, and W.~D. Smith.
\newblock A new table of constant weight codes.
\newblock {\em IEEE Trans. Inform. Theory}, 36(6):1334--1380, 1990.

\bibitem{10}
Y.~M. Chee, S.~H. Dau, A.~C.~H. Ling, and S.~Ling.
\newblock Linear size optimal {$q$}-ary constant-weight codes and constant-composition codes.
\newblock {\em IEEE Trans. Inform. Theory}, 56(1):140--151, 2010.

\bibitem{2}
Y.~M. Chee, F.~Gao, H.~M. Kiah, A.~C.~H. Ling, H.~Zhang, and X.~Zhang.
\newblock Decompositions of edge-colored digraphs: a new technique in the construction of constant-weight codes and related families.
\newblock {\em SIAM J. Discrete Math.}, 33(1):209--229, 2019.

\bibitem{7}
Y.~M. Chee, H.~M. Kiah, and P.~Purkayastha.
\newblock Estimates on the size of symbol weight codes.
\newblock {\em IEEE Trans. Inform. Theory}, 59(1):301--314, 2013.

\bibitem{13}
Y.~M. Chee and S.~Ling.
\newblock Constructions for {$q$}-ary constant-weight codes.
\newblock {\em IEEE Trans. Inform. Theory}, 53(1):135--146, 2007.

\bibitem{14}
W.~Chu, C.~J. Colbourn, and P.~Dukes.
\newblock Constructions for permutation codes in powerline communications.
\newblock {\em Des. Codes Cryptogr.}, 32(1-3):51--64, 2004.

\bibitem{CF07}
D.~J. Costello and G.~D. Forney.
\newblock Channel coding: The road to channel capacity.
\newblock {\em Proceedings of the IEEE}, 95(6):1150--1177, 2007.

\bibitem{18}
T.~Etzion.
\newblock Optimal constant weight codes over {$Z_k$} and generalized designs.
\newblock {\em Discrete Math.}, 169(1-3):55--82, 1997.

\bibitem{5}
T.~Etzion.
\newblock A new approach for examining {$q$}-{S}teiner systems.
\newblock {\em Electron. J. Combin.}, 25(2):Paper No. 2.8, 24, 2018.

\bibitem{12}
G.~Ge.
\newblock Construction of optimal ternary constant weight codes via {B}haskar {R}ao designs.
\newblock {\em Discrete Math.}, 308(13):2704--2708, 2008.

\bibitem{20}
R.~L. Graham and N.~J.~A. Sloane.
\newblock Lower bounds for constant weight codes.
\newblock {\em IEEE Trans. Inform. Theory}, 26(1):37--43, 1980.

\bibitem{J62}
S.~M. Johnson.
\newblock A new upper bound for error-correcting codes.
\newblock {\em IRE Trans.}, IT-8:203--207, 1962.

\bibitem{K96}
J.~Kahn.
\newblock A linear programming perspective on the {F}rankl-{R}\"odl-{P}ippenger theorem.
\newblock {\em Random Structures Algorithms}, 8(2):149--157, 1996.

\bibitem{LS25}
M.~Liu and C.~Shangguan.
\newblock Approximate generalized {S}teiner systems and near-optimal constant weight codes.
\newblock {\em J. Combin. Theory Ser. A}, 209:Paper No. 105955, 19, 2025.

\bibitem{M89}
C.~McDiarmid.
\newblock On the method of bounded differences.
\newblock In {\em Surveys in combinatorics, 1989 ({N}orwich, 1989)}, volume 141 of {\em London Math. Soc. Lecture Note Ser.}, pages 148--188. Cambridge Univ. Press, Cambridge, 1989.

\bibitem{15}
P.~R.~J. \"Osterg\aa~rd and M.~Svanstr\"om.
\newblock Ternary constant weight codes.
\newblock {\em Electron. J. Combin.}, 9(1):Research Paper 41, 23, 2002.

\bibitem{3}
S.~C. Polak.
\newblock Semidefinite programming bounds for constant-weight codes.
\newblock {\em IEEE Trans. Inform. Theory}, 65(1):28--38, 2019.

\bibitem{16}
J.~A. Wood.
\newblock The structure of linear codes of constant weight.
\newblock {\em Trans. Amer. Math. Soc.}, 354(3):1007--1026, 2002.

\bibitem{11}
D.~Wu and P.~Fan.
\newblock Constructions of optimal quaternary constant weight codes via group divisible designs.
\newblock {\em Discrete Math.}, 309(20):6009--6013, 2009.

\bibitem{9}
H.~Zhang and G.~Ge.
\newblock Optimal ternary constant-weight codes of weight four and distance six.
\newblock {\em IEEE Trans. Inform. Theory}, 56(5):2188--2203, 2010.

\bibitem{6}
H.~Zhang and G.~Ge.
\newblock Optimal quaternary constant-weight codes with weight four and distance five.
\newblock {\em IEEE Trans. Inform. Theory}, 59(3):1617--1629, 2013.

\bibitem{8}
H.~Zhang, X.~Zhang, and G.~Ge.
\newblock Optimal ternary constant-weight codes with weight 4 and distance 5.
\newblock {\em IEEE Trans. Inform. Theory}, 58(5):2706--2718, 2012.

\end{thebibliography}

\end{document}